\documentclass{amsart}%
\usepackage{graphicx}
\usepackage{amscd}
\usepackage{amsmath}
\usepackage{amsfonts}
\usepackage{amssymb}%
\providecommand{\U}[1]{\protect\rule{.1in}{.1in}}
\providecommand{\U}[1]{\protect\rule{.1in}{.1in}}
\providecommand{\U}[1]{\protect\rule{.1in}{.1in}}
\newtheorem{theorem}{Theorem}
\theoremstyle{plain}

\newtheorem{corollary}{Corollary}

\newtheorem{definition}{Definition}

\newtheorem{lemma}{Lemma}

\newtheorem{remark}{Remark}

\numberwithin{equation}{section}

\begin{document}
\title{Coincidences for multiple summing mappings}
\author{Geraldo Botelho and Daniel Pellegrino}
\address[G. Botelho]{ Faculdade de Matem\'{a}tica, UFU, Uberl\^{a}ndia, MG, Brazil,
[D.\ Pellegrino] Depto de Matem\'{a}tica, UFPB, Jo\~{a}o Pessoa, PB, Brazil}
\email{botelho@ufu.br, dmpellegrino@gmail.com}
\thanks{A resume of this paper was published in the proceedings of the congress ENAMA
- Second Edition, pages 27-28, Jo\~{a}o Pessoa, PB, Brazil, November 2008}
\maketitle
\maketitle

\section{Notation}

Throughout this paper $n$ is a positive integer, $E_{1},\ldots,E_{n},E$ and
$F$ will stand for Banach spaces over $\mathbb{K}=\mathbb{R}$ or $\mathbb{C}$,
and $E^{\prime}$ is the dual of $E$. By $\mathcal{L}(E_{1},\ldots,E_{n};F)$ we
denote the Banach space of all continuous $n$-linear mappings from
$E_{1}\times\cdots\times E_{n}$ to $F$ with the usual sup norm. If
$E_{1}=\cdots=E_{n}=E$, we write $\mathcal{L}(^{n}E;F)$ and if $F=\mathbb{K}$
we simply write $\mathcal{L}(E_{1},\ldots,E_{n})$ and $\mathcal{L}(^{n}E)$.
Let $p\geq1$. By $\ell_{p}(E)$ we mean the Banach space of all absolutely
$p$-summable sequences $(x_{j})_{j=1}^{\infty}$, $x_{j}\in E$ for all $j$,
with the norm $\Vert(x_{j})_{j=1}^{\infty}\Vert_{p}=\left(  \sum_{j=1}%
^{\infty}\Vert x_{j}\Vert^{p}\right)  ^{1/p}$. $\ell_{p}^{w}(E)$ denotes the
Banach space of all sequences $(x_{j})_{j=1}^{\infty}$, $x_{j}\in E$ for all
$j$, such that $(\varphi(x_{j}))_{j=1}^{\infty}\in\ell_{p}$ for every
$\varphi\in E^{\prime}$ with the norm
\[
\Vert(x_{j})_{j=1}^{\infty}\Vert_{w,p}=\sup\{\Vert(\varphi(x_{j}%
))_{j=1}^{\infty}\Vert_{p}:\varphi\in E^{\prime},\Vert\varphi\Vert\leq1\}.
\]

\begin{definition}
Let $1\leq p_{j}\leq q$, $j=1,\ldots,n$. An $n$-linear mapping $A\in
\mathcal{L}(E_{1},\ldots,E_{n};F)$ is multiple $(q;p_{1},\ldots,p_{n}%
)$-summing if there is a constant $C\geq0$ such that
\[
\left(  \sum_{j_{1},\ldots,j_{n}=1}^{m_{1},\ldots,m_{n}}\left\Vert A(x_{j_{1}%
}^{(1)},\ldots,x_{j_{n}}^{(n)})\right\Vert ^{q}\right)  ^{1/q}\leq
C\Vert(x_{j}^{(1)})_{j=1}^{m_{1}}\Vert_{w,p_{1}}...\Vert(x_{j}^{(n)}%
)_{j=1}^{m_{n}}\Vert_{w,p_{n}}%
\]
for every $m_{1},\ldots,m_{n}\in N$ and any $x_{j_{k}}^{(k)}\in E_{k}$,
$j_{k}=1,\ldots,m_{k}$, $k=1,\ldots,n$. It is clear that we may assume
$m_{1}=\cdots=m_{n}$. The infimum of the constants $C$ working in the
inequality is denoted by $\pi_{q;p_{1},\ldots,p_{n}}(A)$.
\end{definition}

The subspace $\Pi_{q;p_{1},\ldots,p_{n}}^{n}(E_{1},\ldots,E_{n};F)$ of
$\mathcal{L}(E_{1},\ldots,E_{n};F)$ of all multiple $(q;p_{1},\ldots,p_{n}%
)$-summing becomes a Banach space with the norm $\pi_{q;p_{1},\ldots,p_{n}%
}(\cdot)$. If $p_{1}=\cdots=p_{n}=p$ we say that $A$ is multiple
$(q;p)$-summing and write $A\in\Pi_{q;p}^{n}(E_{1},\ldots,E_{n};F)$. The
symbols $\Pi_{q;p_{1},\ldots,p_{n}}^{n}(^{n}E;F)$,$\Pi_{q;p}^{n}(^{n}E;F)$,
\linebreak$\Pi_{q;p_{1},\ldots,p_{n}}^{n}(E_{1},\ldots,E_{n})$, $\Pi_{q;p}%
^{n}(E_{1},\ldots,E_{n})$, $\Pi_{q;p_{1},\ldots,p_{n}}^{n}(^{n}E)$ and
$\Pi_{q;p}^{n}(^{n}E)$ are defined in the obvious way.

Making $n=1$ we recover the classical ideal of absolutely $(q;p)$-summing
linear operators, for which the reader is referred to Diestel, Jarchow and
Tonge \cite{Diestel}. For the space of absolutely $(q;p)$-summing linear
operators from $E$ to $F$ we shall write $\Pi_{q;p}(E;F)$ rather than
$\Pi_{q;p}^{1}(E;F)$.

\section{Results involving cotype $2$}

The next lemma will be used several times.

\begin{lemma}
\label{rr}Let $(a_{j}^{(1)})_{j},...,(a_{j}^{(k)})_{j}\in l_{p}$, $1\leq
p<\infty$. Then $(a_{j_{1}}^{(1)}...a_{j_{k}}^{(k)})_{j_{1},...,j_{k}}\in
l_{p}.$
\end{lemma}

Proof. Let us consider $k=2$. The other cases are similar.
\[%
{\displaystyle\sum\limits_{j_{1},j_{2}=1}^{n}}
\left\vert a_{j_{1}}^{(1)}a_{j_{2}}^{(2)}\right\vert ^{p}=\left(
{\displaystyle\sum\limits_{j_{1}=1}^{n}}
\left\vert a_{j_{1}}^{(1)}\right\vert ^{p}\right)  \left(
{\displaystyle\sum\limits_{j_{2}=1}^{n}}
\left\vert a_{j_{2}}^{(2)}\right\vert ^{p}\right)  .
\]
Making $n\rightarrow\infty$, we conclude that%
\[%
{\displaystyle\sum\limits_{j_{1},j_{2}=1}^{\infty}}
\left\vert a_{j_{1}}^{(1)}a_{j_{2}}^{(2)}\right\vert ^{p}=\left(
{\displaystyle\sum\limits_{j_{1}=1}^{\infty}}
\left\vert a_{j_{1}}^{(1)}\right\vert ^{p}\right)  \left(
{\displaystyle\sum\limits_{j_{2}=1}^{\infty}}
\left\vert a_{j_{2}}^{(2)}\right\vert ^{p}\right)  <\infty.
\]
$\square$\bigskip

The next result is an extension of \cite[Corollary 11.16 (a)]{Diestel}:

\begin{theorem}
\label{a}If $E_{1},...,E_{n}$ have cotype $2$, then%
\[
\Pi_{2}^{n}(E_{1},\ldots,E_{n};F)\subset\Pi_{1}^{n}(E_{1},\ldots,E_{n};F).
\]

\end{theorem}

Proof. Suppose that $A\in\Pi_{2}^{n}(E_{1},\ldots,E_{n};F)$. Let $(x_{j}%
^{(k)})_{j}\in l_{1}^{w}(E_{k}),$ $k=1,...,n.$ Since $E_{k}$ has cotype $2$,
from \cite[Proposition 6]{Oscar},
\[
l_{1}^{w}(E_{k})=l_{2}l_{2}^{w}(E_{k}).
\]
So,
\[
(x_{j}^{(k)})_{j}=(a_{j}^{(k)}y_{j}^{(k)})_{j}\in l_{2}l_{2}^{w}(E_{k}).
\]
Then, using the lemma and H\"{o}lder Inequality,%
\begin{align*}
\left(  A\left(  x_{j_{1}}^{(1)},...,x_{j_{n}}^{(n)}\right)  \right)
_{j_{1},...,j_{n}}  &  =\left(  A\left(  a_{j_{1}}^{(1)}y_{j_{1}}%
^{(1)},...,a_{j_{n}}^{(n)}y_{j_{n}}^{(n)}\right)  \right)  _{j_{1},...,j_{n}%
}\\
&  =\left(  a_{j_{1}}^{(1)}...a_{j_{n}}^{(n)}A\left(  y_{j_{1}}^{(1)}%
,...,y_{j_{n}}^{(n)}\right)  \right)  _{j_{1},...,j_{n}}\in l_{1}\text{
(Holder).}%
\end{align*}
$\square$

A consequence of \cite[Corollary 11.16 (a)]{Diestel} is that if $E$ has cotype
$2,$ then
\[
\Pi_{p}(E;F)=\Pi_{r}(E;F)
\]
for every $1\leq p\leq r\leq2.$ The next theorem shows that for multiple
summing mappings we have a similar situation, except perhaps for the case
$r=2:$

\begin{theorem}
\label{b}If $E_{1},...,E_{n}$ have cotype $2$, then%
\[
\Pi_{p}^{n}(E_{1},\ldots,E_{n};F)=\Pi_{r}^{n}(E_{1},\ldots,E_{n};F)
\]
for every $1\leq p\leq r<2.$
\end{theorem}

Proof. The inclusion $\subset$ is due to David-P\'{e}rez-Garc\'{\i}a
\cite{DavidSTUDIA}.

Suppose that $A\in\Pi_{r}^{n}(E_{1},\ldots,E_{n};F)$. Let $(x_{j}^{(k)}%
)_{j}\in l_{1}^{w}(E_{k})$, $k=1,...,n.$ Since $E_{k}$ has cotype $2$, $E_{k}$
has also cotype $q$ with $r^{\prime}>q>2.$

From \cite[Proposition 6]{Oscar},
\[
l_{1}^{w}(E_{k})=l_{r^{\prime}}l_{r}^{w}(E_{k}).
\]
So,
\[
(x_{j}^{(k)})_{j}=(a_{j}^{(k)}y_{j}^{(k)})_{j}\in l_{r^{\prime}}l_{r}%
^{w}(E_{k}).
\]
Then, using the lemma and H\"{o}lder Inequality,%
\begin{align*}
\left(  A\left(  x_{j_{1}}^{(1)},...,x_{j_{n}}^{(n)}\right)  \right)
_{j_{1},...,j_{n}}  &  =\left(  A\left(  a_{j_{1}}^{(1)}y_{j_{1}}%
^{(1)},...,a_{j_{n}}^{(n)}y_{j_{n}}^{(n)}\right)  \right)  _{j_{1},...,j_{n}%
}\\
&  =\left(  a_{j_{1}}^{(1)}...a_{j_{n}}^{(n)}A\left(  y_{j_{1}}^{(1)}%
,...,y_{j_{n}}^{(n)}\right)  \right)  _{j_{1},...,j_{n}}\in l_{1}\text{
(Holder).}%
\end{align*}
Hence $A\in\Pi_{1}^{n}(E_{1},\ldots,E_{n};F)$ and, from \cite{DavidSTUDIA}, we
have $A\in\Pi_{p}^{n}(E_{1},\ldots,E_{n};F).\square$

\begin{corollary}
\label{yy}If $E_{1},...,E_{n}$ have cotype $2$, then%
\[
\left\{
\begin{array}
[c]{c}%
\Pi_{2}^{n}(E_{1},\ldots,E_{n};F)\subset\Pi_{p}^{n}(E_{1},\ldots
,E_{n};F)\text{ for every }1\leq p\leq2.\\
\Pi_{p}^{n}(E_{1},\ldots,E_{n};F)=\Pi_{r}^{n}(E_{1},\ldots,E_{n};F)\text{ for
every }1\leq p\leq r<2.
\end{array}
\right.
\]

\end{corollary}

\begin{remark}
Note that the previous result extends \cite[Teorema 4.31]{David}.
\end{remark}

\bigskip

\bigskip In \cite{DavidSTUDIA} it is shown that if $1\leq p<q\leq2$ and $F$
has cotype $2$, then%
\[
\Pi_{p}^{n}(E_{1},\ldots,E_{n};F)\subset\Pi_{q}^{n}(E_{1},\ldots,E_{n};F).
\]

Using this result and Corollary \ref{yy}, we have:

\begin{theorem}
If $E_{1},...,E_{n}$ and $F$ have cotype $2$, then%
\[
\Pi_{p}^{n}(E_{1},\ldots,E_{n};F)=\Pi_{r}^{n}(E_{1},\ldots,E_{n};F)
\]
for every $1\leq p\leq r\leq2.$
\end{theorem}

\section{\bigskip Results involving cotype $>2$}

In \cite[Corollary 11.16 (b)]{Diestel} it is shown that if $E$ has cotype
$q>2$, then%
\[
\Pi_{s}(E;F)=\Pi_{r}(E;F)
\]
for every $1\leq s\leq r<\frac{q}{q-1}.$

The next results generalize this result to multiple summing mappings:

\begin{theorem}
If $E_{1},...,E_{n}$ have cotype $q>2$, then%
\[
\Pi_{s}^{n}(E_{1},\ldots,E_{n};F)\subset\Pi_{1}^{n}(E_{1},\ldots,E_{n};F)
\]
for every $1\leq s<\frac{q}{q-1}.$
\end{theorem}

Proof. Since $s<\frac{q}{q-1}$, we have $s^{\prime}>q.$ Suppose that $A\in
\Pi_{s}^{n}(E_{1},\ldots,E_{n};F)$. Let $(x_{j}^{(k)})_{j}\in l_{1}^{w}%
(E_{k})$, $k=1,...,n.$ Since $E_{k}$ has cotype $q>2$, from \cite[Proposition
6]{Oscar},
\[
l_{1}^{w}(E_{k})=l_{s^{\prime}}l_{s}^{w}(E_{k}).
\]
So,
\[
(x_{j}^{(k)})_{j}=(a_{j}^{(k)}y_{j}^{(k)})_{j}\in l_{s^{\prime}}l_{s}%
^{w}(E_{k}).
\]
Then, using Lemma \ref{rr} and H\"{o}lder Inequality, we have%
\begin{align*}
\left(  A\left(  x_{j_{1}}^{(1)},...,x_{j_{n}}^{(n)}\right)  \right)
_{j_{1},...,j_{n}}  &  =\left(  A\left(  a_{j_{1}}^{(1)}y_{j_{1}}%
^{(1)},...,a_{j_{n}}^{(n)}y_{j_{n}}^{(n)}\right)  \right)  _{j_{1},...,j_{n}%
}\\
&  =\left(  a_{j_{1}}^{(1)}...a_{j_{n}}^{(n)}A\left(  y_{j_{1}}^{(1)}%
,...,y_{j_{n}}^{(n)}\right)  \right)  _{j_{1},...,j_{n}}\in l_{1}\text{
(Holder).}%
\end{align*}

$\square$

\begin{corollary}
If $E_{1},...,E_{n}$ have cotype $q>2$, then%
\[
\Pi_{s}^{n}(E_{1},\ldots,E_{n};F)=\Pi_{p}^{n}(E_{1},\ldots,E_{n};F)
\]
for every $1\leq s\leq p<\frac{q}{q-1}.$\bigskip
\end{corollary}

\begin{remark}
Note that there are possible variations of the previous results by exploring
spaces with different cotypes and the Generalized H\"{o}lder Inequality.
\end{remark}

\begin{remark}
The authors have just become aware that the results of this note appear, with
a completely different and independent proof, in a preprint of D. Popa
\cite{Popa}.
\end{remark}

\end{document}